\documentclass[11pt]{amsart}
\usepackage{mathrsfs}
\usepackage{amssymb}

\usepackage{titletoc}
\usepackage{amscd}
\usepackage{amsmath, amssymb}
\usepackage{amsfonts}
\usepackage[colorlinks,linkcolor=blue,citecolor=blue, pdfstartview=FitH]{hyperref}

\numberwithin{equation}{section}

\theoremstyle{plain}
\newtheorem{thm}{Theorem}[section]
\newtheorem{theorem}[thm]{Theorem}
\newtheorem{lemma}[thm]{Lemma}
\newtheorem{corollary}[thm]{Corollary}
\newtheorem{proposition}[thm]{Proposition}

\theoremstyle{definition}

\newtheorem{question}[thm]{Question}

\newtheorem{remark}[thm]{Remark}

\newtheorem{definition}[thm]{Definition}

\newtheorem{example}[thm]{Example}

\newtheorem{defn-thm}[thm]{Definition-Theorem}

\newcommand{\sO}{{\mathcal O}}

\newcommand{\C}{{\mathbb C}}

\renewcommand{\P}{{\mathbb P}}
\newcommand{\Q}{{\mathbb Q}}

\renewcommand{\S}{{\mathbb S}}

\newcommand{\qtq}[1]{\quad\mbox{#1}\quad}
\newcommand{\bp}{\bar{\partial}}
\newcommand{\Om}{\Omega}

\newcommand{\ts}{\otimes}

\newcommand{\btheorem}{\begin{theorem}}
\newcommand{\etheorem}{\end{theorem}}
\newcommand{\bproposition}{\begin{proposition}}
\newcommand{\eproposition}{\end{proposition}}
\newcommand{\bdefinition}{\begin{definition}}
\newcommand{\edefinition}{\end{definition}}
\newcommand{\bcorollary}{\begin{corollary}}
\newcommand{\ecorollary}{\end{corollary}}
\newcommand{\bproof}{\begin{proof}}
\newcommand{\eproof}{\end{proof}}
\newcommand{\bremark}{\begin{remark}}
\newcommand{\eremark}{\end{remark}}
\newcommand{\eexample}{\end{example}}
\newcommand{\bexample}{\begin{example}}
\newcommand{\la}{\langle}
\newcommand{\elemma}{\end{lemma}}
\newcommand{\blemma}{\begin{lemma}}
\newcommand{\ra}{\rangle}
\newcommand{\sq}{\sqrt{-1}}

\newcommand{\p}{\partial}

\renewcommand{\bar}{\overline}
\newcommand{\eps}{\varepsilon}

\renewcommand{\phi}{\varphi}

\newcommand{\ee}{\end{eqnarray*}}
\newcommand{\be}{\begin{eqnarray*}}

\newcommand{\beq}{\begin{equation}}
\newcommand{\eeq}{\end{equation}}

\newcommand{\bd}{\begin{enumerate}}
\newcommand{\ed}{\end{enumerate}}

\renewcommand{\tilde}{\widetilde}

\renewcommand{\>}{\rightarrow}

\usepackage{fancyhdr}
\renewcommand{\leftmark}{{\scriptsize RC-positivity, vanishing theorems and rigidity of holomorphic maps}}

\begin{document}
\title{RC-positivity, vanishing theorems and rigidity of holomorphic maps} \makeatletter
\let\uppercasenonmath\@gobble
\let\MakeUppercase\relax
\let\scshape\relax
\makeatother
\author{Xiaokui Yang}
\date{}
\address{{Address of Xiaokui Yang: Morningside Center of Mathematics, Institute of
        Mathematics, Hua Loo-Keng Center of Mathematical Sciences,
        Academy of Mathematics and Systems Science,
        Chinese Academy of Sciences, Beijing, 100190, China.}}
\email{\href{mailto:xkyang@amss.ac.cn}{{xkyang@amss.ac.cn}}}

\noindent\thanks{This work was partially supported   by China's
Recruitment
 Program of Global Experts and  NSFC 11688101.}
\maketitle

\begin{abstract} Let $M$ and $N$ be two compact complex
manifolds. We show that if the tautological line bundle
$\sO_{T_M^*}(1)$ is not pseudo-effective and $\sO_{T_N^*}(1)$ is
nef, then there is no non-constant holomorphic map from $M$ to $N$.
In particular, we prove that any holomorphic map from a compact
complex manifold $M$ with RC-positive tangent bundle to a compact
complex manifold $N$ with nef cotangent bundle must be a constant
map. As an application, we obtain  that there is no non-constant
holomorphic map from a compact \emph{Hermitian manifold} with
\emph{positive holomorphic sectional curvature} to a Hermitian
manifold with non-positive holomorphic bisectional curvature.

\end{abstract}
\setcounter{tocdepth}{1} \tableofcontents

\section{Introduction}

The classical Schwarz-Pick lemma states that any holomorphic map
from the unit disc in the complex plane into itself decreases the
Poincar\'e metric. This was  extended by Ahlfors (\cite{Ahl38}) to
maps from the disc into a hyperbolic Riemann surface, and by Chern
\cite{Che68} and Lu \cite{Lu68} to higher-dimensional complex
manifolds. A major advance was Yau's Schwarz Lemma \cite{Yau78},
which says that a holomorphic map from a complete K\"ahler manifold
with \emph{Ricci curvature} bounded below into a Hermitian manifold
with holomorphic bisectional curvature bounded above by a negative
constant, is distance decreasing up to a constant depending only on
these bounds. In particular, there is no nontrivial holomorphic map
from compact K\"ahler manifolds with positive Ricci curvature to
 Hermitian manifolds with non-positive holomorphic
bisectional curvature. Later generalizations  were mainly in two
directions: relaxing the curvature hypothesis or the K\"ahler
assumption. In philosophy, holomorphic maps from ``positively
curved" complex manifolds to ``non-positively curved" complex
manifolds should be constant. For more details, we refer to the
recent paper \cite{Tos07} of Tosatti and the references therein.
There are also some other generalizations along this line, for
instance, on complex
analyticity of harmonic maps (e.g. \cite{Siu80, JY93}).\\

In this paper, we obtain a  rigidity theorem on holomorphic maps
between
 complex manifolds, which recovers many classical rigidity theorems along this line in differential geometry. The curvature condition of the domain
manifold is only required to be \emph{RC-positive}. This curvature
notion was introduced in our previous paper \cite{Yang18}, and it is
significantly weaker than the positivity of Ricci curvature. For
instance, a  \emph{complex manifold} with \emph{positive holomorphic
sectional curvature} is RC-positive.
 One of the key ingredients in our proofs relies on the Leray-Grothendieck spectral
sequence and isomorphisms of various cohomology groups, which is
quite different from classical methods in differential geometry. As
it is well-known, the latter is based on various maximum principles
(e.g. \cite{Yau75}).\\

In \cite{Yang18}, we introduced a terminology called
``RC-positivity". A Hermitian holomorphic vector bundle $(\mathscr
E,h^{\mathscr E})$ over a complex manifold $X$ is called
\emph{RC-positive} (resp. RC-negative), if for any $q\in X$ and any
nonzero vector $v\in \mathscr E_q$, there exists \textbf{some}
nonzero  vector $u\in T_qX$ such that $$ R^\mathscr E(u,\bar
u,v,\bar v)>0 \qtq{(resp. $<0$.)}$$ It is easy to see that, for a
Hermitian line bundle $(\mathscr L,h^{\mathscr L})$, it is
RC-positive if and only if its Ricci curvature $-\sq\p\bp\log
h^{\mathscr L}$ has at least one positive eigenvalue at each point
of $X$.
 This terminology has
many nice properties. For instances, quotient bundles of RC-positive
bundles are also RC-positive; subbundles of RC-negative bundles are
still RC-negative. On the other hand, it is obvious that compact
complex manifold with positive holomorphic sectional curvature has
RC-positive tangent bundle. By using Calabi-Yau theorem
\cite{Yau78a}, we proved in \cite[Corollary~3.8]{Yang18} that the
holomorphic tangent bundles of Fano manifolds can admit RC-positive
K\"ahler metrics. This curvature notion should be closely related to
the pseudo-effectiveness of vector bundles defined by P\u{a}un and
Takayama in \cite{PT18} (see also \cite{DPS01}, \cite{Paun16} and
Theorem \ref{equi}). Moreover, it can also be regarded as a
differential geometric interpretation of the positive $\alpha$-slope
investigated by Campana and P\u{a}un in \cite{CP15}. Properties of
RC-positive vector bundles are studied  in \cite{Yang18} and Section
\ref{2}.\\

 The geometry  of vector bundles are usually
 characterized by their tautological line bundles. Let $\mathscr E$ be a
 holomorphic vector bundle and $\P(\mathscr E^*)$ be its projective
 bundle. The tautological line bundle is denoted by $\sO_\mathscr
 E(1)$. For instance, $\mathscr E$ is called ample (resp. nef) if the tautological line
 bundle $\sO_{\mathscr E}(1)$ is ample (resp. nef) over $\P(\mathscr
 E^*)$ (\cite{Har66}). There are many methods to construct Hermitian metrics on  line bundles  (e.g. on $\sO_{\mathscr
 E}(1)$) with various weak positivity. However, it is still a challenge problem to construct
 Hermitian metrics on vector bundles with
 desired curvature properties. For instance, it is a
 long-standing open problem (\cite{Gri69}) to construct positive Hermitian metrics
 on ample vector bundles.\\

 The main result of this paper
is the following rigidity theorem.

\btheorem\label{main1} Let $M$ and $N$ be two compact complex
manifolds. If the tautological line bundle $\sO_{T_M^*}(-1)$ is
RC-positive and $\sO_{T_N^*}(1)$ is  nef.  Then any holomorphic map
from $M$  to $N$ is constant. \etheorem

\noindent Theorem \ref{main1} has an equivalent algebraic version:

\btheorem\label{main0} Let $M$ and $N$ be two compact complex
manifolds. If the tautological line bundle $\sO_{T_M^*}(1)$ is not
pseudo-effective and $\sO_{T_N^*}(1)$ is nef. Then any holomorphic
map from $M$  to $N$ is constant. \etheorem

\noindent Let's explain the curvature conditions in Theorem
\ref{main1} and Theorem \ref{main0}. A line bundle is called
pseudo-effective if it possesses  a (possibly) singular Hermitian
metric whose curvature is semi-positive in the sense of current.
When $M$ is a Riemann surface, $\sO_{T_M^*}(1)$ is not
pseudo-effective if and only if $T^*_M$ is not pseudo-effective,
i.e. $M\cong \P^1$. In this case, Theorem \ref{main0} is classical.
In higher dimensional case, $\sO_{T_M^*}(1)$ is not pseudo-effective
if and only if $\sO_{T_M^*}(1)$ is RC-negative, or equivalently, the
dual line bundle $\sO_{T_M^*}(-1)$ is RC-positive (Theorem
\ref{equi}). Roughly speaking, it says that $T_M$ has a ``positive
direction" at each point of $M$. As we discussed before, the
RC-positivity of $\sO_{T_M^*}(-1)$ is a very weak curvature
condition. For example, it can be implied by  the positivity of
holomorphic sectional curvature (e.g. Proposition \ref{equi}).
Moreover, compact complex manifold with RC positive
$\sO_{T_M^*}(-1)$ is not necessarily K\"ahler. For the curvature
requirement on the target manifold $N$, $\sO_{T_N^*}(1)$ is nef if
and only if the cotangent bundle $T^*_N$ is nef. For instances, all
submanifolds of abelian varieties have nef cotangent bundles.
 The proofs of Theorem \ref{main1} and Theorem \ref{main0} rely on  vanishing theorems for
twisted vector bundles (Theorem \ref{20}) which are establised by
using the Le Potier isomorphism (Leray-Grothendieck spectral
sequence) and characterizations of RC-positive vector bundles
obtained in \cite{Yang17D, Yang18},  which are significantly different from classical methods in differential geometry.\\

 We call that $M$ has \emph{RC-positive tangent bundle} if $M$ admits
 a smooth Hermitian metric $\omega_g$ such that $(T_M,\omega_g)$ is
RC-positive. We show in Proposition \ref{RC-1} that if $T_M$ is
RC-positive, then $\sO_{T_M^*}(-1)$ is RC-positive. As an
application of Theorem \ref{main1}, we obtain the following result.

\btheorem\label{main2} Let $f:M\>N$ be a holomorphic map between two
compact complex manifolds. If $M$ has RC-positive tangent bundle and
$N$ has nef cotangent bundle, then $f$ is a constant map. \etheorem

\noindent There are many  K\"ahler and non-K\"ahler complex
manifolds with RC-positive tangent bundles. We just list some of
them for readers' convenience. \bd
\item[$\bullet$] Fano manifolds \cite[Corollary~3.8]{Yang18};
\item[$\bullet$]   manifolds with positive second Chern-Ricci
curvature \cite[Corollary~3.7]{Yang18};
\item[$\bullet$]  Hopf manifolds $\S^1\times \S^{2n+1}$ (\cite[formula (6.4)]{LY14});

\item[$\bullet$]  complex manifolds with positive holomorphic sectional
curvature. \ed

\noindent The following differential geometric version of Theorem
\ref{main2} is of particular interest, which recovers several
classical rigidity theorems in complex differential geometry.

\bcorollary\label{main3} Let $M$ be a compact complex manifold with
RC-positive tangent bundle, and $N$ be a Hermitian manifold with
non-positive holomorphic bisectional curvature, then any holomorphic
map from $M$ to $N$ is a constant map. \ecorollary

\bremark A compact complex manifold with RC-positive tangent bundle
can  contain no rational curves. For instances,  Hopf manifolds
$\S^1\times \S^{2n+1}$. \eremark

\bremark As we pointed out before,  one of the key ingredients in
the proofs is the Leray-Grothendieck spectral sequence. Although the
conditions in Corollary \ref{main3} are differential geometric,
classical methods (e.g. maximum principles) in differential geometry
can not work for the proof of Corollary \ref{main3}. In Section
\ref{Appendix},
 we  include a discussion on classical methods for
readers' convenience. \eremark

\noindent As a special case of Corollary \ref{main3}, we obtain

\bcorollary\label{main4} Let $(M,\omega_g)$ be a compact Hermitian
manifold with positive holomorphic sectional curvature, and  $(N,
h)$ be a  Hermitian manifold with non-positive holomorphic
bisectional curvature. Then there is no non-constant meropmorphic
map from $M$ or its blowing-up to $N$. \ecorollary

\bremark The notion of positive holomorphic sectional curvature
 is very natural in
differential geometry, but it seems to be mysterious in literature.
Recently, we proved in \cite[Theorem~1.7]{Yang18} that a compact
K\"ahler manifold with positive holomorphic sectional curvature must
be projective and rationally connected, which confirms a well-known
conjecture (\cite[Problem~47]{Yau82}) of S.-T. Yau.  However, the
geometry of \emph{compact complex manifolds} with positive
holomorphic sectional curvature is still not clear.  For some
related topics, we refer to \cite{Ni98, ACH15, Liu16, Yang16, AHZ16,
YZ16, AH17, CY17, Mat18, NZ1, NZ2} and the references therein.  A
project on the geometry of complete non-compact complex manifolds
with RC-positive curvature is also carried out and we have obtained
some results analogous to Yau's classical work \cite{Yau78}.
\eremark

\textbf{Acknowledgements.} I am very grateful to Professor
 K.-F. Liu and Professor S.-T. Yau for their support, encouragement and stimulating
discussions over  years. I would also like to thank Professors J.-P.
Demailly, Y.-X. Li, W.-H. Ou, V. Tosatti, Y.-H. Wu, J. Xiao,
 and X.-Y. Zhou for helpful suggestions.

\vskip 2\baselineskip

\section{Background materials}\label{2}

Let $(\mathscr E,h)$ be a Hermitian holomorphic vector bundle over a
complex manifold $X$ with Chern connection $\nabla$. Let
$\{z^i\}_{i=1}^n$ be the  local holomorphic coordinates
  on $X$ and  $\{e_\alpha\}_{\alpha=1}^r$ be a local frame
 of $\mathscr E$. The curvature tensor $R^{\mathscr E}\in \Gamma(X,\Lambda^{1,1}T^*_X\ts \mathrm{End}(\mathscr E))$ has components \beq R^\mathscr E_{i\bar j\alpha\bar\beta}= -\frac{\p^2
h_{\alpha\bar \beta}}{\p z^i\p\bar z^j}+h^{\gamma\bar
\delta}\frac{\p h_{\alpha \bar \delta}}{\p z^i}\frac{\p
h_{\gamma\bar\beta}}{\p \bar z^j}.\label{cu}\eeq (Here and
henceforth we sometimes adopt the Einstein convention for
summation.) If $(X,\omega_g)$ is a  Hermitian manifold, then
$(T_X,g)$ has Chern curvature components \beq R_{i\bar j k\bar
\ell}=-\frac{\p^2g_{k\bar \ell}}{\p z^i\p\bar z^j}+g^{p\bar
q}\frac{\p g_{k\bar q}}{\p z^i}\frac{\p g_{p\bar \ell}}{\p\bar
z^j}.\eeq The Chern-Ricci curvature $\mathrm{Ric}(\omega_g)$ of
$(X,\omega_g)$ is represented by
$$R_{i\bar j}=g^{k\bar\ell} R_{i\bar j k\bar\ell}$$
and the second Chern-Ricci curvature $\mathrm{Ric}^{(2)}(\omega_g)$
has components
$$R^{(2)}_{k\bar \ell}=g^{i\bar j} R_{i\bar j k\bar\ell}.$$

\bdefinition A Hermitian holomorphic vector bundle $(\mathscr E,h)$
over a complex manifold $X$ is called \emph{Griffiths positive
 at  point $q\in X$}, if for any nonzero vector $v\in \mathscr E_q$, and any nonzero vector $u\in T_qX$
we have  \beq R^\mathscr E(u,\bar u,v,\bar v)>0.\eeq  $(\mathscr
E,h)$ is called \emph{Griffiths positive} if it is Griffiths
positive at every  point of $X$. \edefinition

\noindent As analogous to Griffiths positivity, we introduced in
\cite{Yang18} the following concept.

 \bdefinition\label{Def}
A Hermitian holomorphic vector bundle $(\mathscr E,h)$ over a
complex manifold $X$ is called \emph{RC-positive
 at  point $q\in X$}, if for each nonzero vector $v\in \mathscr E_q$, there exists \textbf{some} nonzero  vector $u\in T_qX$
such that \beq R^\mathscr E(u,\bar u,v,\bar v)>0.\eeq  $(\mathscr
E,h)$ is called \emph{RC-positive} if it is  RC-positive
 at every  point of $X$.  \edefinition

\bremark Similarly, one can define semi-positivity, negativity and etc..
\eremark

\noindent In \cite[Theorem~1.4]{Yang17D}, we obtained an
 equivalent characterization for  RC-positive  line
 bundles which plays a key role in this paper.
 \btheorem\label{equi} Let $X$ be a compact complex manifold and $\mathscr L$ be a holomorphic line bundle over $X$. Then the
 following statements
 are equivalent.

 \bd \item $\mathscr L$ is  RC-positive;
 \item the dual line bundle $\mathscr L^{*}$ is not
 pseudo-effective.
 \ed
 \etheorem

 \noindent As an application of Theorem \ref{equi}, we have

\bcorollary\label{vanline} Let $X$ be a compact complex manifold. If
$\mathscr L$ is an RC-positive line bundle over  $X$, then \beq
H^0(X,\mathscr L^{*})=0.\eeq \ecorollary

\bproof Suppose $H^0(X,\mathscr L^{*})\neq 0$, then $\mathscr L^*$
is $\Q$-effective and so it is pseudo-effective. By Theorem
\ref{equi}, this is a contradiction. \eproof

 The points of the projective bundle $\P(\mathscr E^*)$ of
$\mathscr E\>X$ can be identified with the hyperplanes of $\mathscr
E$. Note that every hyperplane $\mathscr V$ in $\mathscr E_z$
corresponds bijectively to the line of linear forms in $\mathscr
E^*_z$ which vanish on $\mathscr V$.
 Let $\pi : \P(\mathscr E^*) \> X$ be
the natural projection. There is a tautological hyperplane subbundle
$\mathscr  S$ of $\pi^*\mathscr E$ such that
$$\mathscr S_{[\xi]} = \xi^{-1}(0) \subset \mathscr E_z$$ for
all $\xi\in \mathscr E_z^*\setminus\{0\}$. The quotient line bundle
$\pi^*\mathscr E/\mathscr  S$ is denoted $\sO_\mathscr E(1)$ and is
called the \emph{tautological line bundle} associated to $\mathscr
E\>X$. Hence there is an exact sequence of vector bundles over
$\P(\mathscr E^*)$ \beq 0 \>\mathscr S\>\pi^*\mathscr
E\>\sO_\mathscr E(1)\> 0. \label{quotient}\eeq A holomorphic vector
bundle $\mathscr E\>X$ is called \emph{ample} (resp.
\emph{semi-ample, nef}) if the line bundle $\sO_\mathscr
E(1)$ is ample (resp. semi-ample, nef) over $\P(\mathscr E^*)$.\\

 Suppose $\dim_\C X=n$ and
$r=\mathrm{rank}(\mathscr E)$. Let $\pi$ be the projection
$\P(\mathscr E^*)\> X$ and $\mathscr L=\sO_\mathscr E(1)$. Let
$(e_1,\cdots, e_r)$ be the local holomorphic frame on $\mathscr E$
and the dual frame on $\mathscr E^*$ is denoted by $(e^1,\cdots,
e^r)$. The corresponding holomorphic coordinates on $\mathscr E^*$
are denoted by $(W_1,\cdots, W_r)$.  If $(h_{\alpha\bar\beta})$ is
the matrix representation of a smooth Hermitian metric $h^{\mathscr
E}$ on $\mathscr E$ with respect to the basis
$\{e_\alpha\}_{\alpha=1}^r$, then the induced Hermitian metric
$h^\mathscr L$ on $\mathscr L$ can be written as \beq h^\mathscr
L=\frac{1}{\sum h^{\alpha\bar\beta}W_\alpha\bar
W_\beta}\label{inducedmetric} \eeq
 The curvature of $(\mathscr L,h^\mathscr L)$ is \beq
R^{\mathscr L}=\sq\p\bp\log\left(\sum
h^{\alpha\bar\beta}W_\alpha\bar W_\beta\right)
\label{inducedcurvature}\eeq where $\p$ and $\bp$ are operators on
the total space $\P(\mathscr E^*)$. We fix a point $p\in \P(\mathscr
E^*)$, then there exist local holomorphic coordinates
 $(z^1,\cdots, z^n)$ centered at point $q=\pi(p)$ and local holomorphic basis $\{e_1,\cdots, e_r\}$ of $\mathscr E$ around $q$ such that
 \beq h_{\alpha\bar\beta}=\delta_{\alpha\bar\beta}-R^{\mathscr E}_{i\bar j \alpha\bar\beta}z^i\bar z^j+O(|z|^3) \label{normal}\eeq
Without loss of generality, we assume $p$  is the point $(0,\cdots,
0,[a_1,\cdots, a_r])$ with $a_r=1$. On the chart $U=\{W_r=1\}$ of
the fiber $\P^{r-1}$, we set $w^A=W_A$ for $A=1,\cdots, r-1$. By
formula (\ref{inducedcurvature}) and (\ref{normal}) \beq R^{\mathscr
L}(p)=\sq\sum R^{\mathscr E}_{i\bar j\alpha\bar\beta}\frac{a_\beta
\bar a_\alpha}{|a|^2}dz^i\wedge d\bar z^j+\omega_{\mathrm{FS}}
\label{induced} \eeq where
$|a|^2=\sum\limits_{\alpha=1}^r|a_\alpha|^2$ and
$\omega_{\mathrm{FS}}=\sq
\sum\limits_{A,B=1}^{r-1}\left(\frac{\delta_{AB}}{|a|^2}-\frac{a_B\bar
a_A}{|a|^4}\right)dw^A\wedge d\bar w^B$ is the Fubini-Study metric
on the fiber $\P^{r-1}$. The following result is one of the key
ingredients in this paper.

\bproposition\label{RC-1} Let $X$ be a compact complex manifold. If
$(\mathscr E, h^{\mathscr E})$ is an RC-positive vector bundle over
$X$, then $\sO_{\mathscr E^*}(-1)$ is an RC-positive line bundle
over $\P(\mathscr E)$. \eproposition  \bproof By (\ref{induced}),
the induced metric on $\sO_{\mathscr E^*}(-1)$ over $\P(\mathscr E)$
has curvature form
$$ R^{\sO_{\mathscr E^*}(-1)}=-\left(\sq\sum
R^{\mathscr E^*}_{i\bar j\alpha\bar\beta}\frac{a_\beta \bar
a_\alpha}{|a|^2}dz^i\wedge d\bar z^j+\omega_{\mathrm{FS}}\right). $$
On the other hand, $R^{\mathscr E^*}=-\left(R^{\mathscr E}\right)^t$
and so $$ R^{\sO_{\mathscr E^*}(-1)}=\sq\sum R^{\mathscr E}_{i\bar
j\alpha\bar\beta}\frac{a_\alpha \bar a_\beta}{|a|^2}dz^i\wedge d\bar
z^j-\omega_{\mathrm{FS}}. $$ Hence,
 $\sO_{\mathscr
E^*}(-1)$ is RC-positive    as long as $(\mathscr E, h^{\mathscr
E})$ is RC-positive . \eproof

\bremark We also have the following results.

\bd \item If $\mathscr L_1$ is an RC-positive line bundle and
$\mathscr L_2$ is a pseudo-effective line bundle, then $\mathscr
L_1\ts \mathscr L_2$ is RC-positive;

\item Let $(\mathscr E,h^{\mathscr E})$ be an RC-positive vector bundle and
$(\mathscr F,h^{\mathscr F})$ be a Griffiths semi-positive vector
bundle. The  Hermitian vector bundle $(\mathscr E\ts \mathscr F,
h^{\mathscr E}\ts h^{\mathscr F})$ is \emph{not} necessarily
RC-positive unless $\mathrm{rank}(\mathscr{E})=1$. \ed

\eremark

\bremark It is easy to see that if $(\mathscr E, h^{\mathscr E})$ is
Griffiths positive (resp. semi-positive), then the tautological line
bundle $\sO_{\mathscr E}(1)$ is positive (resp. semi-positive). It
is a long-standing open problem (so called Griffiths conjecture)
that whether the converse is valid. In the same vein, we wonder
whether the RC-positivity of $\sO_{\mathscr E^*}(-1)$ can imply the
RC-positivity of $\mathscr E$. \eremark

The following well-known lemma is called \emph{the Le Potier
isomorphism}(\cite{LeP75}). Its proof relies on the
Leray-Grothendieck spectral sequence, and we refer to
\cite[Theorem~5.16]{SS85} and the references therein.

\blemma\label{Leiso} Let $\mathscr E$ be a holomorphic vector bundle
over a complex manifold $X$ and $\mathscr F$ be a coherent sheaf on
$X$. Then for all $p,q\geq 0$ \beq H^q(X,\Om_X^p\ts \mathscr E\ts
\mathscr F)\cong H^q\left(\P(\mathscr E^*),\Om_{\P(\mathscr
E^*)}^p\ts \sO_{\mathscr E}(1)\ts \pi^*\mathscr F\right),\eeq where
$\pi:\P(\mathscr E^*)\>X$ is the projection. In particular, \beq
H^0(X,\mathscr E)\cong H^0(\P(\mathscr E^*), \sO_{\mathscr
E}(1)).\label{Liso}\eeq
 \elemma

\noindent By using the Le Potier isomorphism, we obtain vanishing
theorems for vector bundles. \blemma\label{vanishing} Let $\mathscr
E$ be a holomorphic vector bundle over  a compact complex manifold
$X$. If $\sO_{\mathscr E^*}(-1)$ is RC-positive, then \beq
H^0(X,\mathscr E^*)=0. \label{van}\eeq In particular, if $\mathscr
E$ is RC-positive, then $\mathscr E^*$ has no nontrivial holomorphic
section. \elemma \bproof It follows from Corollary \ref{vanline},
Proposition \ref{RC-1} and Lemma \ref{Leiso}. \eproof

\noindent The following concept is a generalization of the
RC-positivity for line bundles.
\begin{definition}\label{Def22} Let $\mathscr L$ be a holomorphic line bundle
over a  complex manifold $X$.  $\mathscr L$ is called
\emph{$k$-positive}, if there exists a smooth Hermitian metric
$h^{\mathscr L}$  such that the Chern curvature $R^{{\mathscr
L}}=-\sq\p\bp\log h^{\mathscr L}$ has at least $(\dim X-k)$ positive
eigenvalues at every point on $X$.
\end{definition}
It is easy to see that $\mathscr L$ is $(\dim X-1)$-positive if and
only if it is RC-positive.  In \cite[Theorem~14]{AG62}, Andreotti
and Grauert proved the following fundamental vanishing theorem.
\blemma \label{AG} Let $\mathscr L$ be a $k$-positive line bundle
over a compact complex manifold $X$. Then for any vector bundle
$\mathscr F$ over $X$, there exists a positive integer
$m_0=m_0(\mathscr F)$ such that \beq H^q(X,\mathscr L^{\ts m}\ts
\mathscr F)=0\eeq for all $q>k$ and $m\geq m_0$. \elemma

\vskip 2\baselineskip

\section{Vanishing theorems for tensor product of vector bundles}

\noindent The main result of this section is the following vanishing
theorem.

\btheorem\label{20} Let $\mathscr E$ and $\mathscr F$ be two
holomorphic vector bundles over a compact complex manifold $X$. If
$\sO_{\mathscr E^*}(-1)$ is RC-positive over $\P(\mathscr E)$ and
$\sO_{\mathscr F}(1)$ is nef over $\P(\mathscr F^*)$. Then \beq
H^0(X,\mathscr E^*\ts \mathscr F^*)=0. \label{va20} \eeq \etheorem

\noindent By Theorem \ref{equi}, we have a variant of Theorem
\ref{20}. \btheorem\label{200} Let $\mathscr E$ and $\mathscr F$ be
two holomorphic vector bundles over a compact complex manifold $X$.
If $\sO_{\mathscr E^*}(1)$ is not pseudo-effective  and
$\sO_{\mathscr F}(1)$ is nef. Then \beq H^0(X,\mathscr E^*\ts
\mathscr F^*)=0. \eeq \etheorem

\bremark  Theorem \ref{20} does not hold in general if
$\sO_{\mathscr F}(1)$ is pseudoeffective and $\mathrm{rank}(\mathscr
F)>1$. It should hold if we refine this notion a bit more (e.g.
\cite{PT18,DPS01}). \eremark

\noindent By using Proposition \ref{RC-1}, we have
\btheorem\label{2000} Let $\mathscr E$ and $\mathscr F$ be two
holomorphic vector bundles over a compact complex manifold $X$. If
$\mathscr E$ is RC-positive and $\mathscr F$ is nef, then \beq
H^0(X,\mathscr E^*\ts \mathscr F^*)=0.  \eeq \etheorem

\noindent Before giving the proof of Theorem \ref{20}, we need
several  lemmas.

\blemma\label{pbrc} Let $f:X\>Y$ be a holomorphic \emph{submersion}
between two complex manifolds. If $\mathscr L$ is an RC-positive
line bundle over $Y$, then $f^*(\mathscr L)$ is also RC-positive.
\elemma \bproof Suppose $\dim X=m$ and $\dim Y=n$. Let
$\{z^i\}_{i=1}^m$ and $\{w^\alpha\}_{\alpha=1}^n$ be the local
holomorphic coordinates on $X$ and $Y$ respectively. Let $h$ be a
smooth RC-positive metric on $\mathscr L$ and $R=-\sq\p\bp\log h$.
It is easy to see that the curvature tensor of $(f^*(\mathscr L),
f^*h)$ is given by \beq R_{\alpha\bar\beta} \frac{\p f^\alpha}{\p
z^i}\frac{\p \bar f^\beta}{\p\bar z^j}dz^i\wedge d\bar
z^j.\label{7}\eeq Since $(\mathscr L, h)$ is RC-positive, at any
point $p\in Y$, there exists a nonzero local vector $v=(v^1,\cdots,
v^n)$ such that $\sum R_{\alpha\bar\beta}v^\alpha\bar v^\beta>0$.
Since $f$ is a smooth submersion, the rank of the matrix
$\left(\frac{\p f^\alpha}{\p z^i}\right)$ is equal to $n=\dim Y$.
Therefore, there exists a nonzero vector $u=(u^1,\cdots u^m)$ such
that $\left(\frac{\p f^\alpha}{\p z^i}\right)u=v.$  Hence,
$(f^*(\mathscr L), f^*h)$ is RC-positive.
 \eproof
\bremark Lemma \ref{pbrc} also holds for $k$-positive line bundles.
\eremark

\blemma\label{pbsp} Let $f:X\>Y$ be a holomorphic map between two
compact complex manifolds. If $\mathscr L$ is a nef line bundle over
$Y$, then $f^*(\mathscr L)$ is also nef. \elemma \bproof It follows
from the definition of nefness and formula (\ref{7}). \eproof

\blemma\label{pbspv} Let $f:X\>Y$ be a holomorphic map between two
compact complex manifolds. If  $\mathscr E$ is a holomorphic vector
bundle over $Y$ such that $\sO_{\mathscr E}(1)$ is nef. Then
$\sO_{f^* \mathscr E}(1)$ is also nef. \elemma \bproof We have the
following commutative diagram
$$\CD
  \sO_{f^* \mathscr E}(1) @>(f^\#)^*>> \sO_{\mathscr E}(1) \\
  @V  VV @V  VV  \\
  \P(f^*(\mathscr E^*)) @>f^\#>> \P(\mathscr E^*) \\
  @V  VV @V  VV  \\
  X @>f>> Y.
\endCD
$$
Lemma \ref{pbspv} follows from the above diagram and Lemma
\ref{pbsp}.\eproof

\blemma \label{555} Let $\mathscr{E}$ be a holomorphic vector bundle
over a compact complex manifold $X$. If  $\sO_{\mathscr{E}}(1)$ is
$(\dim X-1)$-positive over $\P(\mathscr{E}^*)$, then \beq
H^0(X,\mathscr E^*)=0 \eeq \elemma

 \bproof  If
$\sO_{\mathscr{E}}(1)$ is  $(\dim X-1)$-positive over
$\P(\mathscr{E}^*)$, then by Lemma \ref{AG}, for any vector bundle
$\mathscr F$ on $\P(\mathscr E^*)$, there exists some positive
integer $m_0=m_0(\mathscr F)$ such that
$$H^q(\P(\mathscr E^*), \sO_{\mathscr E}(m)\ts \mathscr F)=0$$
for all $q>\dim X-1$ and $m\geq m_0$. In particular, if we take
$q=n=\dim X$ and $\mathscr F=\Om^n_{\P(\mathscr E^*)}$, by Lemma
\ref{Leiso} and the Serre duality,
$$H^n(\P(\mathscr E^*), \sO_{\mathscr E}(m)\ts \Om^n_{\P(\mathscr E^*)})\cong H^n(X,\mathrm{Sym}^{\ts m}\mathscr E\ts \Om_X^n)\cong H^0(X, \mathrm{Sym}^{\ts m}\mathscr E^*)=0.$$
In particular, for large $m$, we have
$$H^0(\P(\mathscr E), \sO_{\mathscr E^*}(m))=0.$$
Hence, $H^0(\P(\mathscr E), \sO_{\mathscr E^*}(1))=0$ and
$H^0(X,\mathscr E^*)=0$.
 \eproof

\btheorem\label{key} Let $X$ be a compact complex manifold.  If
$(\mathscr{L},h^{\mathscr L})$ is an RC-positive line bundle, and
$\mathscr{E}$ is a holomorphic vector bundle with nef tautological
line bundle $\sO_{\mathscr{E}}(1)$. Then
$$H^0\left(X,\mathscr E^*\ts \mathscr{L}^*\right)=0.$$ \etheorem

\bproof Let $ \pi:\P(\mathscr{E^*})\>X$ be the natural projection.
Since $\pi$ is a submersion, by Lemma \ref{pbrc},
$\pi^*\mathscr{L}$ is RC-positive.\\

 \emph{Claim $1$.} $\pi^*\mathscr{L}\ts \sO_{\mathscr{E}}(1)$ is a
$(\dim X-1)$-positive line bundle over $\P(\mathscr{E^*})$.\\

\noindent  Fix a smooth Hermitian metric $h^{\mathscr{E}}$ on
$\mathscr{E}$ and a smooth Hermitian metric $\omega$ on
$\P(\mathscr{E^*})$. The induced metric on $\sO_{\mathscr{E}}(1)$ is
denoted by $h^{\sO_{\mathscr{E}}(1)}$. Since the restriction of
$h^{\sO_{\mathscr{E}}(1)}$ on each fiber $\P^{r-1}$ is a
Fubini-Study metric, by curvature formula (\ref{induced}), there
exist a Hermitian metric $\omega_X$ on $X$ and two positive
constants $c_1, c_2$ such that \beq -\sq\p\bp\log h^{\sO_{\mathscr
E}(1)}+c_1\pi^*(\omega_X)\geq c_2\omega. \label{3.6} \eeq

Let $\lambda(x)$ be the largest eigenvalue function of the curvature
tensor $-\sq\p\bp\log h^{\mathscr L}$ of $(\mathscr{L}, h^{\mathscr
L})$ with respect to the Hermitian metric $\omega_X$ on $X$ and \beq
c_3=\min_{x\in X} \lambda(x). \eeq Since $X$ is compact and
$-\sq\p\bp\log h^{\mathscr L}$ is RC-positive, we deduce $c_3>0$.
Moreover, at any point $q\in X$, there exists a nonzero vector
$u_0\in T_q X$ such that \beq \left(-\sq\p\bp\log h^{\mathscr
L}\right)(u_0,\bar u_0)\geq c_3|u_0|^2_{\omega_X}. \eeq

\noindent Since $\pi^*:\P(\mathscr E^*)\>X$ is a holomorphic
submersion, by Lemma \ref{pbrc}, $h_1=\pi^*(h^{\mathscr{L}})$ is an
RC-positive  metric on $\pi^*\mathscr{L}$. Moreover, for any point
$p\in\P(\mathscr E^*)$ with $\pi(p)=q\in X$, there exists a nonzero
vector $u_1\in T_p\P(\mathscr E^*)$ such that $\pi_*(u_1)=u_0\in T_q
X$ and \beq  \left(-\sq\p\bp\log h_1\right)(u_1,\bar u_1)=
\left(-\sq\p\bp\log h^{\mathscr L}\right)(u_0,\bar u_0)\geq
c_3|u_0|^2_{\omega_X}>0. \label{1} \eeq
 We fix a small number $\eps>0$
such that \beq \frac{c_3}{2}-c_1\eps>0. \label{3}\eeq

\noindent On the other hand, since $\sO_{\mathscr E}(1)$ is nef,
there exists a smooth Hermitian metric $h_0$ on
$\sO_{\mathscr{E}}(1)$ such that the curvature of $(\sO_{\mathscr
E}(1),h_0)$ satisfies \beq -\sq\p\bp\log h_0\geq -\eps c_2\omega
\label{3.11} \eeq
 over
$\P(\mathscr{E^*})$. Let
$h=\left(h^{\sO_{\mathscr{E}}(1)}\right)^\eps\cdot h_0^{(1-\eps)}$
be a smooth Hermitian metric on $\sO_{\mathscr{E}}(1)$. Then
$$\left(\pi^*\mathscr{L}\ts \sO_{\mathscr{E}}(1),
h_1\ts h\right)$$ is  $(\dim X-1)$-positive, i.e. the curvature
tensor $R=-\sq\p\bp\log (h_1h)$ has at least $r$-positive
eigenvalues at each point of $\P(\mathscr{E^*})$. Indeed, \beq
R=\eps\left(-\sq\p\bp\log
h^{\sO_{\mathscr{E}}(1)}\right)+(1-\eps)\left(-\sq\p\bp\log
h_0\right)+ \left(-\sq\p\bp\log h_1\right).\label{3.12} \eeq
 By (\ref{3.6}), (\ref{1}), (\ref{3}),
(\ref{3.11}) and (\ref{3.12})  we have \be R(u_1,u_1)&\geq&
\eps(c_2|u_1|^2_{\omega}-c_1|u_1|^2_{\pi^*\omega_X})-(1-\eps)\eps
c_2 |u_1|^2_\omega +c_3|u_1|^2_{\pi^*\omega_X} \\&\geq&
\frac{c_3}{2}|u_1|^2_{\pi^*\omega_X}=\frac{c_3}{2}|u_0|^2_{\omega_X}
>0. \ee Along the fiber $\P^{r-1}$ direction, for any $u_2\in T_p
\P(\mathscr E^*)$ with $\pi_*(u_2)=0\in T_qX$, we have
$$R(u_2,u_2)\geq
\eps(c_2|u_2|^2_{\omega}-c_1|u_2|^2_{\pi^*\omega_X})-(1-\eps)\eps
c_2 |u_2|^2_\omega+c_3|u_2|^2_{\pi^*\omega_X}\geq c_2\eps^2
|u_2|^2_\omega.$$ Since the map $\pi_*:T_p\P(\mathscr E)\>T_q X$ is
surjective,  $\dim \mathrm{ker}(\pi_*)=r-1$ and $u_1\notin
\mathrm{ker}(\pi_*)$, we deduce that the curvature tensor
$R=-\sq\p\bp\log (h_1h)$ has at least $r$ positive eigenvalues at
each point of $\P(\mathscr{E^*})$.\\

\emph{Claim $2$.} The tautological line bundle $\sO_{\mathscr L\ts
\mathscr E}(1)$ is $(\dim X-1)$-positive over $\P(\mathscr L^*\ts
\mathscr E^*)$. Indeed, it follows from the fact that $i:\P(\mathscr
L^*\ts \mathscr E^*)\>\P(\mathscr{E^*})$ is an isomorphism and \beq
\sO_{\mathscr L\ts \mathscr E}(1)=i^*\left(\sO_{\mathscr{E}}(1)\ts
\pi^*(\mathscr L) \right). \eeq

 \noindent By Lemma \ref{555}, we obtain $H^0\left(X,\mathscr E^*\ts
 \mathscr{L}^*\right)=0$. The proof of Theorem \ref{key} is
 completed.
 \eproof

\noindent\emph{The proof of Theorem \ref{20}.} Let
$\pi:\P(\mathscr{E})\>X$ be the  projection. By Lemma \ref{Leiso}
\beq H^0(X,\mathscr E^*\ts \mathscr F^*)\cong H^0\left(X,
\pi_*\left(\sO_{\mathscr E^*}(1)\ts \pi^* \mathscr
F^*\right)\right)\cong H^0\left(Y,\sO_{\mathscr E^*}(1)\ts \pi^*
\mathscr F^*\right), \label{4}\eeq where $Y=\P(\mathscr E)$. Since
$\sO_{\mathscr F}(1)$ is nef, by Lemma \ref{pbspv}
$\sO_{\pi^*\mathscr F}(1)$ is nef. Let $\mathscr L=\sO_{\mathscr
E^*}(-1)$, $ \mathscr W={\pi^*\mathscr F}$ and $\tilde \pi:
\P(\mathscr W)\>Y$. Since $\mathscr L$ is an RC-positive line bundle
and $\sO_{\mathscr W}(1)$ is nef,  by   Lemma \ref{Leiso} and
Theorem \ref{key}, \beq H^0(X,\mathscr F^*\ts \mathscr E^*)\cong
H^0\left(Y,{\pi^*\mathscr F^*}\ts \sO_{\mathscr E^*}(1)\right)=
H^0(Y,\mathscr W^*\ts \mathscr L^*)=0.\label{6}\eeq The proof of
Theorem \ref{20} is completed.\qed

\vskip 2\baselineskip

\section{RC-positivity and rigidity of holomorphic maps }

In this section, we prove the main results of this paper, i.e.
Theorem \ref{main1} (=Theorem \ref{main55}), Theorem
\ref{main2}(=Theorem \ref{main50}) and Corollary \ref{main4}
(=Corollary \ref{main60}).

\btheorem\label{main55} Let $M$ and $N$ be two compact complex
manifolds.  If $\sO_{T_M^*}(-1)$ is an RC-positive line bundle and
$\sO_{T_N^*}(1)$ is  nef.  Then any holomorphic map from $M$  to $N$
is constant. \etheorem

\bproof Let $\mathscr E=T_M\otimes f^*(T^*_N)$ and $\{z_i\}$,
$\{w_\alpha\}$ be the local holomorphic coordinates on $M$ and $N$,
respectively. Let
$$ s = \p f=f^{\alpha}_i dz_i\otimes e_\alpha \in \Gamma(M,\mathscr E^*)$$
where $e_{\alpha} = f^{\ast }\frac{\partial }{\partial w_{\alpha}}$.
Since $f$ is a holomorphic map, $s$ is a holomorphic section of
$\mathscr E$, i.e. $s\in H^0(M,\mathscr E^*)$. Since
$\sO_{T^*_N}(1)$ is nef, by Lemma \ref{pbspv}, we know
$\sO_{f^*(T^*_N)}(1)$ is also nef. By Theorem \ref{20},
$H^0(M,\mathscr E^*)=0$. Hence $s$ is a constant map. \eproof

\noindent In particular, we have

\btheorem\label{main50} Let $M$ be a compact complex manifold with
RC-positive tangent bundle $T_M$ and $N$ be a compact complex
manifold with nef cotangent bundle. Then any holomorphic map from
$M$ to $N$ is constant. \etheorem \bproof By Proposition \ref{RC-1},
if $T_M$ is RC-positive, then $\sO_{T_M^*}(-1)$ is RC-positive.
Theorem \ref{main50} follows from Theorem \ref{main55}. \eproof

Let $M, N$ be compact complex manifolds of complex dimensions $m$
and $n$ respectively. Recall that a \emph{meromorphic map} $f:M\>N$
is given by an irreducible analytic subset (the graph of $f$)
$\Gamma\subset M\times N$ together with a proper analytic subset
$S\subset M$ and a \emph{holomorphic} map $f:M\setminus S\>N$ such
that $\Gamma$ restricted to $(M-S)\times N$ is exactly the graph of
$f$.

\bcorollary\label{main60} Let $(M,\omega_g)$ be a compact Hermitian
manifold with positive holomorphic sectional curvature, and  $(N,
h)$ be a  Hermitian manifold with non-positive holomorphic
bisectional curvature. Then there is no non-constant meropmorphic
map from $M$ or its blowing-up to $N$. \ecorollary

\bproof Let $f:M\>N$ be a meromorphic map. By a theorem of Griffiths
(\cite[Theorem~II]{Gri71}) and Shiffman (\cite[Theorem~2]{Shi71}),
when the target manifold has non-positive holomorphic sectional
curvature, then $f$ is  holomorphic.  It is easy to see that if
$\omega_g$ has positive holomorphic sectional curvature, then
$(T_M,\omega_g)$ is RC-positive. By Theorem \ref{main50}, there is
no non-constant holomorphic map from $M$ to $N$.

Let $\tilde M$ be a blowing-up of $M$ along some submanifold and
$\pi:\tilde
 M\>M$ be the canonical map. If $\tilde f:\tilde M\>N$ is a meromorphic
 map, then it is holomorphic. Moreover, it induces a meromorphic map
 $f:M\>N$. Hence $f$ is a constant map. By Aronszajin's principle (\cite{Aro57}),  $\tilde
 f$ is also constant.
 \eproof

\noindent We have shown in \cite[Corollary~3.7]{Yang18} that if a
complex manifold has positive second Chern-Ricci curvature, then it
is RC-positive.

\bcorollary\label{main5000} Let $f:M\>N$ be a holomorphic map
between two compact complex manifolds. If $(M, g)$ has positive
second Chern-Ricci curvature $\mathrm{Ric}^{(2)}(g)$ and $N$ has nef
cotangent bundle, then $f$ is a constant. \ecorollary

\noindent Hence,  the following classical result is a special case
of Corollary \ref{main5000} (e.g. \cite{YZ18}).

\bcorollary\label{classical} Let  $(M, g)$ be a compact Hermitian
manifold with positive second Chern-Ricci curvature
$\mathrm{Ric}^{(2)}(g)$ and and  $(N, h)$ be a Hermitian manifold
with non-positive holomorphic bisectional curvature, then $f$ is a
constant. \ecorollary

\noindent We would like to propose questions on rigidity results for
more general target manifolds. For instance,
\begin{question} Let $f:M\>N$ be a holomorphic map
between two compact complex manifolds. If $T_M$ is RC-positive (or
$\sO_{T_M^*}(-1)$ is RC-positive) and $N$ is Kobayashi hyperbolic
(or more generally, $N$ is a complex manifold without rational
curves), is $f$ necessarily a constant map?

\end{question}

\vskip 2\baselineskip

\section{Appendix: Yau's Schwarz calculation and rigidity of holomorphic
maps}\label{Appendix}

In this section, we review classical differential geometric methods
(a model version of Yau's Schwarz calculation) on the proof of
rigidity of holomorphic maps. We shall see clearly that the main
results  in this paper (e.g. Corollary \ref{main3}) can not be
proved by using purely differential geometric methods. The following
result is essentially well-known (e.g. \cite{Che68, Lu68, Yau78}).

\begin{lemma}\label{Lemma4.1}
Let $f:(M,g)\rightarrow (N,h)$ be a holomorphic map between two
Hermitian manifolds. Then in the local holomorphic coordinates
$\{z^i\}$ and $\{w^\alpha\}$ on $M$ and $N$, respectively, we have
the identity
$$
\p\bp u = \la \nabla df, \nabla df\ra + \left( R^g_{i\bar j k\bar
\ell} g^{k\bar q} g^{p\bar \ell} h_{\alpha\bar\beta} f^{\alpha }_p
\overline{f^{\beta}_q}  -
R^h_{\alpha\bar\beta\gamma\bar\delta}\left(
 f^{\alpha}_i \overline{f^{\beta}_j} \right)  \left(g^{p\bar q}  f^{\gamma}_p \overline{
f^{\delta}_q} \right)\right)dz^i\wedge d\bar z^j. $$ and $$ \Delta_g
u = |\nabla df|^2 + \left(g^{i\bar j} R^g_{i\bar j k\bar \ell}
\right) g^{k\bar q} g^{p\bar \ell} h_{\alpha\bar\beta} f^{\alpha }_p
\overline{f^{\beta}_q}  -
R^h_{\alpha\bar\beta\gamma\bar\delta}\left(g^{i\bar j}
 f^{\alpha}_i \overline{f^{\beta}_j} \right)  \left(g^{p\bar q}  f^{\gamma}_p \overline{
f^{\delta}_q} \right),$$
 where $u=\mathrm{tr}_{\omega_g} (f^*\omega_h ) $, $f^{\alpha
}_i = \frac{\partial f^{\alpha } } {\partial z^i}$, where $f$ is
represented by $w^{\alpha }=f^{\alpha }(z)$ locally, $\nabla$ is the
induced connection on the bundle $\mathscr E=T^*_M\otimes f^*(T_N)$.
\end{lemma}

\noindent To simplify the formulations, at a given point $p\in M$
and $q=f(p)\in N$, we choose $g_{i\bar j}=\delta_{ij}$ and
$h_{\alpha\bar\beta}=\delta_{\alpha\bar\beta}$. Hence, we have \beq
\p\bp u = \la \nabla df, \nabla df\ra +\left(
\sum_{k,\ell,\alpha}R^g_{i\bar j k\bar \ell}
 f^{\alpha }_k
\overline{f^{\alpha}_\ell}  - \sum_{\alpha,\beta,\gamma,\delta, k}
R^h_{\alpha\bar\beta\gamma\bar\delta}\left(
 f^{\alpha}_i \overline{f^{\beta}_j} \right)  \left(  f^{\gamma}_k \overline{
f^{\delta}_k} \right)\right)dz^i\wedge d\bar z^j\label{new}\eeq and

\beq \Delta_g u = |\nabla df|^2 + \sum_{k,\ell,\alpha}R^{(2)}_{k\bar
\ell}
 f^{\alpha }_k
\overline{f^{\alpha}_\ell}  - \sum_{\alpha,\beta,\gamma,\delta, k,i}
R^h_{\alpha\bar\beta\gamma\bar\delta}\left(
 f^{\alpha}_i \overline{f^{\beta}_i} \right)  \left(  f^{\gamma}_k \overline{
f^{\delta}_k} \right)\label{Yau} \eeq where $R^{(2)}_{k\bar
\ell}=g^{i\bar j} R^g_{i\bar j k\bar\ell}$. If $M$ is compact, by
applying the standard maximum principle to (\ref{Yau}), we obtain
Corollary \ref{classical}. One may wonder whether Corollary
\ref{main3} can be obtained by applying a similar maximum principle
to  equation (\ref{new}). Suppose $u$ attains a maximum at some
point $p\in X$. Then for any vector $v=(v^1,\cdots v^n)$, by formula
(\ref{new}), at point $p\in X$, we have \beq 0\geq \sum_{i,j}
\frac{\p^2 u}{\p z^i\p\bar z^j} v^i\bar v^j\geq
\sum_{i,j}\left(\sum_{k,\ell,\alpha}R^g_{i\bar j k\bar \ell}
 f^{\alpha }_k
\overline{f^{\alpha}_\ell}\right) v^i\bar v^j. \label{diff}\eeq
Recall that, if $(T_M, g)$ is RC-positive, then for any nonzero
vector $\xi=(\xi^1,\cdots, \xi^n)$, there exists some nonzero vector
$\eta=(\eta^1,\cdots, \eta^n)$ (\emph{it may depend on $\xi$}!) such
that $R^g_{i\bar j k\bar\ell} \eta^i\bar \eta^j \xi^k\bar
\xi^\ell>0.$ Apparently, in (\ref{diff}), there are many vectors
indexed by $\alpha$, and in general there does not exist a uniform
vector $v$ such that the right hand side of (\ref{diff}) is
positive. A refined notion called ``uniform RC-positivity"  would
work for this analytical proof. By using similar ideas, we also
investigated rigidity of harmonic maps into Riemannian
manifolds in \cite{Yang18c}. \\

 The relationship  between the Leray-Grothendieck
spectral sequence in algebraic geometry and maximum principles in
differential geometry will be systematically investigated.

\end{document}